\documentclass{amsart}
\usepackage{amssymb, latexsym}

\DeclareMathOperator{\im}{Im}

\DeclareMathOperator{\re}{Re}

\newcommand{\eps}{{\varepsilon}}
\newcommand{\D}{{\mathbb D}}

\newcommand{\C}{{\mathbb C}}
\newcommand{\R}{{\mathbb R}}
\newcommand{\N}{{\mathbb N}}

\newcommand{\hoh}{{\dot{H}^{1/2}}}
\newcommand{\even}{{(e)}}
\newcommand{\odd}{{(o)}}
\newcommand{\plus}{{(+)}}
\newcommand{\minus}{{(-)}}
\newcommand{\plusminus}{{(\pm)}}
\newcommand{\any}{{(\ast)}}
\newcommand{\Szego}{{Szeg\H{o} }}
\newcommand{\interval}{{[-2,2]}}
\newcommand{\W}{{\mathcal W}}

\theoremstyle{plain}
\newtheorem{theorem}{Theorem}
\newtheorem{lemma}[theorem]{Lemma}

\newtheorem{proposition}[theorem]{Proposition}

\numberwithin{equation}{section} \numberwithin{theorem}{section}

\begin{document}

\title[A strong Szeg\H{o} theorem for Jacobi matrices]
{A strong Szeg\H{o} theorem for Jacobi matrices}
\author{E. ~Ryckman}
\address{University of California, Los Angeles CA 90095, USA}
\email{eryckman@math.ucla.edu}

\begin{abstract}
We use a classical result of Gollinski and Ibragimov to prove an
analog of the strong \Szego theorem for Jacobi matrices on
$l^2(\N)$. In particular, we consider the class of Jacobi matrices
with conditionally summable parameter sequences and find necessary
and sufficient conditions on the spectral measure such that
$\sum_{k=n}^\infty b_k$ and $\sum_{k=n}^\infty (a_k^2 - 1)$ lie in
$l^2_1$, the linearly-weighted $l^2$ space.
\end{abstract}

\maketitle

%
%
%
%

\section{Introduction}
Let us begin with some notation.  We study the spectral theory of
Jacobi matrices, that is semi-infinite tridiagonal matrices
$$J = \begin{pmatrix}
b_1 & a_1 & 0 & 0 \\
a_1 & b_2 & a_2 & 0 \vphantom{\ddots}\\
0 & a_2 & b_3 & \ddots\\
0 & 0 & \ddots & \ddots
\end{pmatrix}$$
where $a_n > 0$ and $b_n \in \R$.  In this paper we make the
overarching assumption that the sequences $b_n$ and $a_n^2-1$ are
conditionally summable. We may then define
\begin{align}\label{lambda kappa def}
\begin{split}
\lambda_n &:= -\sum_{k=n+1}^\infty b_k\\
\kappa_n &:= -\sum_{k=n+1}^\infty (a_k^2 - 1)
\end{split}
\end{align}
for $n = 0,1,\dots$.

Let $d\nu$ be the spectral measure for the pair $(J,\delta_1)$,
where $\delta_1 = (1,0,0,\dots)^t$, and assume that $d\nu$ is not
supported on a finite set of points (we will call such measures
\emph{nontrivial}). Let
\begin{equation}\label{m-function}
m(z) := \langle \delta_1 , (J-z)^{-1}\delta_1 \rangle = \int
\frac{d\nu(x)}{x-z}
\end{equation}
be the associated $m$-function, defined for $z \in \C \backslash
\text{supp}(\nu)$.

We will write
$$\{\beta_n\} \in l^2_{s} \quad\text{if}\quad \|\beta\|_{l^2_{s}}^2 := \sum_n
|n|^s|\beta_n|^2 < \infty,$$ and let $\hoh (\partial \D)$ denote
the (homogeneous) Sobolev space of order $1/2$ of functions
defined on $\partial \D$:
$$f \in \hoh \quad\text{if}\quad \|f\|_{\hoh}^2 := \|\hat f (n)\|_{l^2_1}^2
= \sum_n |n| |\hat f (n)|^2 < \infty.$$ If $f$ is a function on
$\interval$, we say $f \in \hoh$ if $f(2 \cos \theta) \in \hoh
(\partial \D).$  Also, we will say $v \in \W$ if $v(x)$ is
supported in $\interval$ and has one of the forms
\begin{equation}\label{weight forms}
\Bigl(\sqrt{4-x^2}\Bigr)^{\pm 1} v_0 (x) \quad \text{or}\quad
\Biggl(\sqrt{\frac{2-x}{2+x}}\Biggr)^{\pm 1} v_0 (x)
\end{equation}
with $\log(v_0) \in \hoh$.

Our main result is:

\begin{theorem}\label{main theorem}
Let $J$ be a Jacobi matrix.  The following are equivalent:
\begin{enumerate}
\item The sequences associated to $J$ by \eqref{lambda kappa def}
obey $\lambda,\kappa \in l^2_1$

\item $J$ has finitely-many eigenvalues that all lie in $\R
\setminus [-2,2]$, and on $\interval$ the spectral measure is
purely absolutely continuous, $d\nu(x) = v(x) dx$, with $v \in
\W$.
\end{enumerate}
\end{theorem}

The main ingredient in the proof will be the following version of
the strong \Szego theorem\footnote{The version we use is due to
\cite{Golinskii Ibragimov} and \cite{Ibragimov}.  For relevant
definitions see, for instance, \cite{Simon}.}.

\begin{theorem}[Golinskii-Ibragimov]\label{lemma Golinskii-Ibragimov}
Let $d\mu$ be a probability measure on $\partial \D$ that is not
supported on a finite set of points, and let $\{\alpha_n\}
\subseteq \D$ be the associated Verblunsky coefficients. The
following are equivalent:
\begin{enumerate}
\item $\alpha \in l^2_1$

\item $d\mu = w \frac{d\theta}{2\pi}$ and $\log w \in \hoh$.
\end{enumerate}
\end{theorem}

We now outline the proof of Theorem \ref{main theorem}. To apply
the strong \Szego theorem we must move to the circle, so we must
first remove all the eigenvalues in $\R \setminus \interval$. To
do so we use double commutation (see \cite{Gesztesy Teschl}):

\begin{theorem}[Double Commutation]\label{lemma double commutation}
Let $E \in \R \setminus \sigma (J)$, and let $\gamma
> 0$. Define a new Jacobi matrix $\widetilde J$ by
\begin{gather*}
\tilde a_n = a_n \frac{\sqrt{c_{n-1}c_{n+1}}}{c_n}\\
\tilde b_n = b_n + \gamma
\Bigl(\frac{a_{n-1}\phi_{n-1}\phi_n}{c_{n-1}}-\frac{a_n \phi_n
\phi_{n+1}}{c_n}\Bigr)
\end{gather*}
where $J \phi = E \phi$, $\phi_0 = 0$, $\phi_1 = 1$ and $$c_n = 1
+ \gamma \sum_{j=1}^n |\phi_j|^2.$$ Then $\sigma (\widetilde J ) =
\sigma (J) \cup \{E\}$, $E$ is a simple eigenvalue of $\widetilde
J$, and $$\tilde m(z) = \frac{1}{1+\gamma}\Bigl(m(z) -
\frac{\gamma}{z-E}\Bigr).$$

Conversely, let $|E|>2$ be a simple eigenvalue of $J$ with
eigenvector $\phi$. Choose $\gamma = -1/\|\phi\|^2$ and define a
new Jacobi matrix $\widetilde J$ as above.  Then $\sigma
(\widetilde J) = \sigma (J) \setminus \{E\}$ and
$$\widetilde m(z) = (1 + \gamma)m(z) + \frac{\gamma}{z-E}.$$
\end{theorem}

We prove an asymptotic integration result in Section 2, which we
combine with the above theorem in Section 3 to prove

\begin{proposition}\label{lemma remove eits}
Let $J$ be a Jacobi matrix, and let $E$ be an isolated eigenvalue
of $J$ in $\R \setminus \interval$.  Let $\widetilde J$ be the
Jacobi matrix obtained from $J$ by removing the eigenvalue $E$
using Theorem \ref{lemma double commutation}.  Then
\begin{enumerate}
\item $\tilde \lambda, \tilde \kappa \in l^2_1$ if and only if
$\lambda, \kappa \in l^2_1$

\item $\tilde v \in \W$ if and only if $v \in \W$.
\end{enumerate}
\end{proposition}

This proposition essentially allows us to consider Theorem
\ref{main theorem} under the additional hypothesis $\sigma (J)
\subseteq \interval$.  This allows us to move to the circle, as
follows. Given a nontrivial probability measure $d\mu$ on
$\partial \D$ that is invariant under complex conjugation, one can
define a nontrivial probability measure $d\nu$ on $[-2,2]$ by
$$\int_{-2}^2 g(x) d\nu(x) = \int_0^{2\pi} g(2 \cos \theta) d\mu
(\theta).$$ Similarly, given such a measure $d\nu$, one can define
a measure $d\mu$ that is symmetric under complex conjugation by
$$\int_0^{2\pi} h(\theta) d\mu(\theta) = \int_{-2}^2 h(\arccos
(x/2)) d\nu(x)$$ when $h(-\theta) = h(\theta)$.

The map $d\mu \mapsto d\nu$ is one of a family of four maps that
we call the \Szego mappings\footnote{The map $Sz^\even$ is due to
\cite{Szego}, while the other three are due to \cite{spaniards},
then developed further in \cite{Killip Nenciu} and \cite{Simon}.}.
We denote it by  $d\nu = Sz^\even(d\mu)$. The other three maps are
given by
\begin{align}\label{szego mappings}
\begin{split}
Sz^\odd (d\mu) &= c^2 (4-x^2) Sz^\even(d\mu)\\
Sz^\plusminus (d\mu) &= c_\pm^2 (2 \mp x) Sz^\even(d\mu)
\end{split}
\end{align}
\begin{gather}\label{c definition}
\begin{split}
c = \frac{1}{\sqrt{2(1-|\alpha_0|^2)(1 - \alpha_1)}}\\
c_\pm = \frac{1}{\sqrt{2(1 \mp \alpha_0)}}.
\end{split}
\end{gather}
If $d\mu$ is absolutely continuous with respect to Lebesgue
measure we will write $d\mu(\theta) = w(\theta)
\frac{d\theta}{2\pi}$ and then $Sz^{\any}(d\mu)(x) = v^\any(x)
dx$. In this case the above relations become
\begin{align}\label{weight relation}
\begin{split}
v^\even(x) &= \frac{1}{\pi \sqrt{4-x^2}}w(\arccos(x/2))\\
v^\odd(x) &= \frac{c}{\pi} \sqrt{4-x^2}w(\arccos(x/2))\\
v^\plusminus(x) &= c_\pm \sqrt{\frac{2\mp x}{2 \pm
x}}w(\arccos(x/2)).
\end{split}
\end{align}

For $\ast \in \{e,o,+,-\}$, we will write $J^\any$ for the Jacobi
matrix determined by $d\nu^\any$ and $a^\any$, $b^\any$ for its
parameter sequences. The relationship between $\alpha$ and
$a^\any,b^\any$ is given by

\begin{proposition}[Direct Geronimus Relations\footnote{The relationship
between $\alpha$ and $a^\even,b^\even$ was first discovered by
\cite{Geronimus}.  The other three were later found by
\cite{spaniards} using techniques similar to \cite{Geronimus}.
\cite{Killip Nenciu} and \cite{Simon} have a different proof using
operator techniques.}]\label{lemma Geronimus relations} Let $d\mu$
be a nontrivial probability measure on $\partial \D$ that is
invariant under conjugation, and let $d\nu^\any = Sz^\any (d\mu)$.
Then for all $n \geq 0$
\begin{gather*}
[a^\even_{n+1}]^2 = (1-\alpha_{2n-1})(1-\alpha_{2n}^2)(1+\alpha_{2n+1})\\
b^\even_{n+1} = \alpha_{2n}(1-\alpha_{2n-1}) -
\alpha_{2n-2}(1+\alpha_{2n-1})
\end{gather*}
\begin{gather*}
[a^\odd_{n+1}]^2 = (1+\alpha_{2n+1})(1-\alpha_{2n+2}^2)(1-\alpha_{2n+3})\\
b^\odd_{n+1} = -\alpha_{2n+2}(1+\alpha_{2n+1})+
\alpha_{2n}(1-\alpha_{2n+1})
\end{gather*}
\begin{gather*}
[a^\plusminus_{n+1}]^2 = (1\pm\alpha_{2n})(1-\alpha_{2n+1}^2)(1\mp\alpha_{2n+2})\\
b^\plusminus_{n+1} = \mp\alpha_{2n+1}(1\pm\alpha_{2n}) \pm
\alpha_{2n-1}(1\mp\alpha_{2n}).
\end{gather*}
\end{proposition}

Since $a_n > 0$, there is no ambiguity in which sign to choose for
the square root above.  We always take $\alpha_{-1} = -1$. The
value of $\alpha_{-2}$ is irrelevant since it is multiplied by
zero.

From the Direct Geronimus Relations we see that decay of the
$\alpha$'s determines decay of the $a$'s and $b$'s. This allows us
to prove one direction of Theorem \ref{main theorem} in Section 4.

To prove the other direction, we will find certain relationships
between the Verblunsky parameters and solutions of $J u = E u$ at
$E = \pm 2$. We study asymptotics of these solutions in Section 5,
then find the desired relationships in Section 6, which we term
the Inverse Geronimus Relations.  In Section 7 we review some Weyl
theory, and in Section 8 we combine all these ideas to finish the
proof.

It is a pleasure to thank Rowan Killip for his helpful advice.

%
%
%
%

\section{Asymptotic integration}
Suppose $\widetilde J$ and $J$ are related through double
commutation (as in Theorem \ref{lemma double commutation}). In the
next section we will relate $\tilde \lambda, \tilde \kappa$ to
$\lambda, \kappa$. By Theorem \ref{lemma double commutation} we
see\footnote{In order to avoid excessive subscripting later in
this section, we will write $a(n)$ for $a_n$, etc.}
\begin{gather*}
|\tilde \kappa(n-1) - \kappa(n-1)| = \Biggl| \sum_{k=n}^\infty
 a(k)^2 \Biggl(\frac{c(k-1)c(k+1)}{c(k)^2}-1\Biggr)\Biggr|\\
|\tilde \lambda(n-1) - \lambda(n-1) | =  |\gamma| \Biggl|
\sum_{k=n}^\infty
\Biggl(\frac{a(k-1)\phi(k-1)\phi(k)}{c(k-1)}-\frac{a(k) \psi(k)
\psi(k+1)}{c(k)}\Biggr)\Biggr|.
\end{gather*}
So to prove part $(1)$ of Proposition \ref{lemma remove eits}, we
must determine asymptotics for $\phi$ when $E \in \R \setminus
\interval$.  To do so we use the theory of asymptotic integration
as developed in \cite{Harris Lutz, Harris Lutz unified, Hartman
Wintner, Levinson} and particularly \cite{Benzaid Lutz}. However,
as we need $l^p_s$ control of the errors (rather than the usual
$o(1)$ control) we must modify their results. Throughout, we will
use the notation $x \lesssim y$ if there is a constant $c
> 0$ such that $x \leq cy$. Also, if $x_n$ is a sequence, we write
$x = y + l^p_s$ to indicate $x_n = y_n + \eps_n$ for some other
sequence $\eps \in l^p_s$.

\begin{proposition}\label{lemma discrete HW off}
Let $\Lambda(k) = diag[\lambda_1(k),\dots,\lambda_n(k)]$ and
suppose that there exists $0 < \delta < 1$ so that for a fixed $i$
either
\begin{equation}\label{dichotomy condition off}
(I)\text{  } \Biggl| \frac{\lambda_i(k)}{\lambda_j(k)} \Biggr|
\geq 1+\delta \quad\text{or}\quad (II)\text{  } \Biggl|
\frac{\lambda_i(k)}{\lambda_j(k)} \Biggr| \leq 1-\delta
\end{equation}
for each $j\neq i$, where $k \geq k_0$ for some $k_0$. Suppose
also that $\|V(k)\| \in l^2_s$ for some $s \geq 0$. Then the
system
\begin{equation}\label{diag system off}
\Psi(k+1) = [ \Lambda(k) + V(k) ] \Psi(k)
\end{equation}
has a solution of the form
\begin{equation}\label{diag solution off}
\Psi_i(k) = \Biggl(\prod_{l = k_0}^{k-1} \lambda_i (l) +
V_{ii}(l)\Biggr) (e_i + l^2_s)
\end{equation}
where $e_i$ is the $i^{th}$ standard unit vector in $\R^n$.
\end{proposition}
As all norms on a finite dimensional space are equivalent, it does
not matter which we mean when we write things like $\|V(k)\| \in
l^p_s$ or $(e_i + l^p_s)$.

We will prove Proposition \ref{lemma discrete HW off} by using a
Harris-Lutz transformation followed by a Levinson-type result.  We
will state and use these results, then prove them at the end of
this section.

\begin{proposition}\label{lemma discrete HL off}
With the assumptions of Proposition \ref{lemma discrete HW off},
there exists a sequence of matrices $Q(k)$ such that $Q(k)_{ii}
=0$, $\|Q(k)\| \in l^2_s$, and
\begin{equation}\label{HL equality off}
V(k) - diag V(k) + \Lambda(k)Q(k) - Q(k+1)\Lambda(k) = 0.
\end{equation}
\end{proposition}

\begin{proposition}\label{lemma discrete L off}
Say $\Lambda(k)$ satisfies the assumptions of Proposition
\ref{lemma discrete HW off}, and suppose that $\|R(k)\| \in l^1_s$
for some $s \geq 0$.  Then the system
\begin{equation}\label{diag system L off}
x(k+1) = [ \Lambda(k) + R(k) ] x(k)
\end{equation}
has a solution of the form
\begin{equation}\label{diag solution L off}
x_i(k) = \Biggl(\prod_{l = k_0}^{k-1} \lambda_i (l) \Biggr) (e_i +
l^2_s).
\end{equation}
\end{proposition}

Assuming Propositions \ref{lemma discrete HL off} and \ref{lemma
discrete L off} we have

\begin{proof}[Proof of Proposition \ref{lemma discrete HW off}]
Let $Q(k)$ be as guaranteed by Proposition \ref{lemma discrete HL
off}, and define $x(k)$ by $$\Psi(k) = [I + Q(k)]x(k)$$ (as $Q(k)
\rightarrow 0$, $[I + Q(k)]$ is invertible for large $k$, so the
above definition makes sense).  Then $\Psi$ is a solution of
\eqref{diag system off} if and only if $x$ solves $$x(k+1) =
[\widetilde \Lambda (k) + \widetilde V (k)] x(k)$$ where
\begin{gather*}
\widetilde \Lambda (k) = \Lambda (k) + diag V(k)\\
\widetilde V(k) = [I + Q(k)]^{-1}[V(k)Q(k) - Q(k+1)diag V(k)].
\end{gather*}
It is easy to see that $\widetilde \Lambda$ still satisfies the
dichotomy condition \eqref{dichotomy condition off}.
 Moreover, as $\|V(k)\| , \|Q(k)\| \in l^2_s$ we have that $\|\widetilde
V(k)\| \in l^{1}_s$. So we may apply Proposition \ref{lemma
discrete L off} to the $x$-system to find a solution
$$x_i(k) = \Biggl(\prod_{l=k_0}^{k-1} \tilde \lambda_i(l)
\Biggr)(e_i + \eps(k))$$ for some $\eps(k) \in l^2_s$. But then
$$\Psi_i(k) = [I + Q(k)]x_i(k) = \Biggl(\prod_{l = k_0}^{k-1}
\lambda_i (l)  + V(l)_{ii}\Biggr) \Bigl(e_i + \eps(k) + Q(k)e_i +
Q(k)\eps(k) \Bigr).$$  By Proposition \ref{lemma discrete HL off}
we have that $\|\eps(k) + Q(k)e_i + Q(k)\eps(k)\| \in l^2_s$, as
required.
\end{proof}

Next we prove Propositions \ref{lemma discrete HL off} and
\ref{lemma discrete L off}.  In doing so we will make frequent use
of the following two lemmas.

\begin{lemma}\label{lemma sum products}
Let $s \geq 1$, $\beta, \gamma \in l^2_{s}$, and define a sequence
$\eta_n := \sum_{k=n}^\infty \beta_k \gamma_k$.  Then $\eta \in
l^2_{s}$ and $\|\eta\|_{l^2_{s}} \leq
\|\beta\|_{l^2_{s}}\|\gamma\|_{l^2_{s}}$.  In particular, if $\tau
\in l^1_s$ then $\sum_{k=n}^\infty \tau_k \in l^2_s$.
\end{lemma}

\begin{proof}
Throughout the proof, all norms refer to $l^2_s$.  By
Cauchy-Schwarz we have
\begin{align*} \|\eta\|^2 &= \sum_{n=1}^\infty
n^{s}\Bigl|\sum_{k=n}^\infty \beta_k \gamma_k \Bigr|^2 \leq
\sum_{n=1}^\infty n^{s} \Bigl(\sum_{k=n}^\infty |\beta_k
\gamma_k|\Bigr)^2 \\
&\leq \sum_{n=1}^\infty n^{s} \Bigl(\sum_{k=n}^\infty |\beta_k|^2
\Bigr) \Bigl(\sum_{k=n}^\infty |\gamma_k|^2 \Bigr) =
\sum_{n=1}^\infty \Bigl(\sum_{k=n}^\infty n^{s} |\beta_k|^2
\Bigr) \Bigl(\sum_{k=n}^\infty |\gamma_k|^2 \Bigr) \\
&\leq \|\beta\|^2 \sum_{n=1}^\infty \sum_{k=n}^\infty |\gamma_k|^2
= \|\beta\|^2 \sum_{k=1}^\infty k |\gamma_k|^2 \leq \|\beta\|^2
\|\gamma\|^2.
\end{align*}
The last statement follows by applying the above argument to
$\beta = \gamma = |\tau|^{1/2}$.
\end{proof}

\begin{lemma}\label{lemma convolution type}
Suppose $\gamma(k) \in l^2_1$ and $\beta > 1$.  Then $
\sum_{l=1}^{k-1} \beta^{2(l-k)}\gamma(l)  \in l^2_1$.
\end{lemma}

\begin{proof}
We must show that $$\gamma \mapsto  \sum_{l=1}^{k-1}
\beta^{2(l-k)}\gamma(l)$$ maps $l^2_1 \rightarrow l^2_1$.
Equivalently we will show $$\gamma \mapsto \sum_{l=1}^{k-1}
\sqrt{\frac{k}{l}}\beta^{2(l-k)}\gamma(l)$$ maps $l^2 \rightarrow
l^2$.  This is an integral operator with kernal
$$h(l,k) = \chi_{\{1,\dots,k-1\}}(l) \sqrt{\frac{k}{l}}
\beta^{2(l-k)}$$ so by Schur's Test this will be a bounded
operator if we can show $$\sup_k \sum_{l=1}^\infty h(l,k) \leq C
\quad\text{and}\quad \sup_l \sum_{k=1}^\infty h(l,k) \leq C$$ for
some $C \geq 0$.  This is done by the following lemma.
\end{proof}

\begin{lemma}
For any $\alpha \in \R$ and $\eps > 0$ we have $$\sup_l
\sum_{k=1}^\infty \Bigl( \frac{|k|+1}{|l|+1}\Bigr)^\alpha
e^{-\eps|k-l|}<\infty.$$
\end{lemma}

The proof is standard and proceeds by splitting the sum at $k=l$
and bounding each piece separately.  We omit the details.

\begin{proof}[Proof of Proposition \ref{lemma discrete HL off}]
Define $Q(k)$ by $Q(k)_{ii} = 0$ and
\begin{gather*}
Q(k)_{ij} = -\sum_{m=k}^\infty \frac{V(m)_{ij}}{\lambda_j(m)}
\prod_{l=k}^m \frac{\lambda_j(l)}{\lambda_i(l)}
\quad\text{if}\quad
(i,j)\in(I)\\
Q(k)_{ij} = \sum_{m=k_0}^{k-1} \frac{V(m)_{ij}}{\lambda_i(m)}
\prod_{l=k}^m \frac{\lambda_i(l)}{\lambda_j(l)}
\quad\text{if}\quad (i,j)\in(II).
\end{gather*}
As $\|V(k)\| \in l^p_s \subseteq l^\infty$, $Q(k)_{ij}$ is
dominated (in either case above) by a convergent geometric series,
so the sums defining $Q$ converge.  By the above definition,
\eqref{HL equality off} holds.

To show that $\|Q(k)\| \in l^p_s$ we argue as follows.  For
$(i,j)\in(I)$ we have that $$\Biggl|
\frac{1}{\lambda_j(m)}\prod_{l=k}^m
\frac{\lambda_i(l)}{\lambda_j(l)}\Biggr| \in l^2_1.$$ Similarly,
for $(i,j)\in(II)$ we have that
$$\Biggl| \frac{1}{\lambda_i(m)}\prod_{l=m}^{k-1}
\frac{\lambda_i(l)}{\lambda_j(l)}\Biggr| \lesssim |\beta|^{m-k}$$
for some $|\beta|>1$. So by Lemmas \ref{lemma sum products} and
\ref{lemma convolution type} we see that $Q \in l^p_s$.
\end{proof}

\begin{proof}[Proof of Proposition \ref{lemma discrete L off}]
Define $w(k)$ by $$x(k) = \Bigl( \prod_{l=k_0}^{k-1} \lambda_i(l)
\Bigr) w(k).$$  Then $x$ solves \eqref{diag system L off} if and
only if $w$ solves the system $$w(k+1) =
\frac{1}{\lambda_i(k)}[\Lambda(k) + R(k)]w(k).$$  We'll compare
the $w$-system to the diagonal system $$y(k+1) =
\frac{1}{\lambda_i(k)}\Lambda(k)y(k).$$ The $y$-system has a
fundamental matrix $$Y(k) = diag
\Biggl[\prod_{l=k_0}^{k-1}\frac{\lambda_1(l)}{\lambda_i(l)},
\dots,1,\dots,\prod_{l=k_0}^{k-1}\frac{\lambda_n(l)}{\lambda_i(l)}\Biggr]$$
with a $1$ in the $i^{th}$ spot. Let $P_1 = diag [p_1,\dots,p_n]$
where $$p_j = \begin{cases} 1,
\quad (i,j)\in(I)\\
0,\quad (i,j) \in (II)
\end{cases}$$ and let $P_2 = I - P_1$.  By the assumptions on
$\Lambda(k)$ we see that for $k_0 \leq l \leq k-1$, $\|Y(k)P_1
Y(l+1)^{-1}\| \leq C$ and for $k_0 \leq k \leq l$, $\|Y(k)P_2
Y(l+1)^{-1}\| \leq C$ for some $C > 0$.

Now let $k_1 \geq k_0$ to be chosen later, and consider the
operator $$[Tz](k) =
\sum_{l=k_1}^{k-1}Y(k)P_1Y(l+1)^{-1}\frac{1}{\lambda_i(l)}R(l)z(l)
- \sum_{l=k}^\infty
Y(k)P_2Y(l+1)^{-1}\frac{1}{\lambda_i(l)}R(l)z(l)$$ acting on
$l^\infty (\N,\C^n)$.  Choose $k_1$ so that
$$2C\delta^{-1}\sum_{l=k_1}^{\infty}\|R(l)\| < \eps < 1$$ (which is
possible because $\|R(k)\| \in l^1_s$).  Then we see that
$$\|Tz\|_{l^\infty} \leq \Biggl( 2 C \delta^{-1}
\sum_{l=k_1}^\infty \|R(l)\|\Biggr) \|z\|_{l^\infty} \leq \eps
\|z\|_{l^\infty}$$ for all $z \in l^\infty$.  Thus, $T : l^\infty
\rightarrow l^\infty$ is a contraction.  In particular, given $y
\in l^\infty$, there exists a unique $w \in l^\infty$ solving $w =
y + Tw$.

Say $y\in l^\infty$ and $w = y + Tw$.  By the definition of $T$,
$y$ is a solution of the $y$-system if and only if $w$ is a
solution to the $w$-system. In particular this holds for $y =
e_i$. It remains to show $w = y + l^2_1$, for which we consider
each of the sums defining $Tw$ separately.  As $\|R(k)\| \in l^1_s
\subseteq l^2_s$ and $\|Y(k)P_1
Y(l+1)^{-1}\frac{1}{\lambda_i(l)}w(l)\| \lesssim 1$, Lemma
\ref{lemma convolution type} shows the first sum is in $l^2_1$.
Similarly, because $\|R(k)\| \in l^1_1$, Lemma \ref{lemma sum
products} shows that the second sum is in $l^2_1$.
\end{proof}

Finally, we allow perturbed diagonalizable systems, rather than
just the perturbed diagonal systems of Proposition \ref{lemma
discrete HW off}.

\begin{proposition}\label{lemma discrete HW off final}
Suppose $A(k)$ has eigenvalues $\lambda_i(k)$ satisfying
\eqref{dichotomy condition off} and $\sup_k |\lambda_j(k)| \leq C$
for all $j$. Let $A(k) = S(k)^{-1}\Lambda(k)S(k)$ where
$\Lambda(k) = diag[\lambda_1(k), \dots,\lambda_n(k)]$, and suppose
that $S(k) \rightarrow S(\infty)$ where $S(\infty)$ is invertible
and $\|S(k+1)-S(k)\| \in l^2_s$ for some $s \geq 0$.  Finally,
suppose $V \in l^2_s$. Then the system
\begin{equation}\label{diag system off final}
\Psi(k+1) = [A(k) + V(k)]\Psi(k)
\end{equation}
has a solution of the form $$\Psi_i(k) =
S(k)^{-1}\Biggl(\prod_{l=k_0}^{k-1} \lambda_i(l) + \widetilde
V(l)_{ii} \Biggr) (e_i + l^2_s)$$ where $$\widetilde V(k) = S(k)
V(k) S(k)^{-1} + \bigl( S(k+1) - S(k) \bigr) \bigl( A(k) + V(k)
\bigr) S(k)^{-1}$$ so in particular $\| \widetilde V (k)\| \in
l^2_s$.
\end{proposition}

\begin{proof}
We'll reduce to the case of Proposition \ref{lemma discrete HW
off}.  Define $z(k) = S(k) \Psi(k)$, so $\Psi$ is a solution of
\eqref{diag system off final} if and only if $z$ solves the system
\begin{equation}\label{z solution off}
z(k+1) = [\Lambda(k) + \widetilde V(k)] z(k)
\end{equation}
where $\widetilde V$ is as in the statement of the proposition.
Now
\begin{align*}
\|\widetilde V(k)\| &\lesssim \|V(k)\| + \|S(k+1) - S(k)\| \bigl(
\|A(k)\| + \|V(k)\|\bigr)\\
&\lesssim \|V(k)\| + \|S(k+1) - S(k)\| \in l^2_s
\end{align*}
because $S(\infty)$ is invertible and $\sup_{j,k} |\lambda_j(k)|
\leq C$.  So by Proposition \ref{lemma discrete HW off}, there
exists a solution to \eqref{z solution off} of the form $$z_i(k) =
\Biggl( \prod_{l=k_0}^{k-1} \lambda_i(l) + \widetilde V(l)_{ii}
\Biggr) (e_i + l^2_s).$$ Undoing the transformation we find a
solution to \eqref{diag system off final} of the form $$\Psi_i(k)
= S(k)^{-1}z_i(k)$$ as desired.
\end{proof}

%
%
%
%

\section{The double commutation result}
In this section we prove Proposition \ref{lemma remove eits}.

\begin{proof}[Proof of Proposition \ref{lemma remove eits}(2)]
Suppose that $\widetilde J$ and $J$ are related by double
commutation at $E \in \R$ an isolated point of $\sigma (J)$. Write
$d\nu(x) = v(x)dx$ and recall that Lebesgue almost everywhere
$$v(x) = \frac{1}{\pi}\im m(x + i0) .$$ By Theorem \ref{lemma
double commutation}
$$\tilde m(z) = \frac{1}{1+\gamma}\Bigl(m(z) -
\frac{\gamma}{z-E}\Bigr).$$ But then
$$\tilde v (x) = \frac{1}{\pi} \im \tilde m(x+i0) =
\frac{1}{\pi(1 + \gamma)} \im  m(x + i0)  = \frac{1}{1+\gamma}
v(x)$$ almost everywhere.  Clearly $\tilde v \in \W$ if and only
if $v \in \W$.
\end{proof}

Part (1) is more difficult, and will take the rest of this section
to prove.  We will use the asymptotic integration results obtained
in Section 2.

\begin{lemma}\label{lemma eit solutions}
Write $E = \beta + \beta^{-1}$ with $|\beta| > 1$.  The recurrence
equation at $E$ has solutions of the form $$\psi_\pm (k) = c_\pm
\beta^{\pm k} (1 + l^2_1)$$ for some constants $c_\pm \in
\R\setminus \{0\}$.
\end{lemma}

\begin{proof}
We will prove the result for $E > 2$, the proof for $E < -2$ being
similar. We can write the recurrence equation
$$a(k+1)\psi(k+1) + \bigl( b(k) - E\bigr)\psi(k) + a(k) \psi (k-1)
= 0$$ as the system
\begin{equation}\label{eigenvalue system}
\Psi(k+1) = [ A(k) + V(k) ]\Psi(k)
\end{equation}
where $$\Psi(k) =
\begin{bmatrix}
\psi(k)\\
\psi(k-1)
\end{bmatrix}
\quad A(k) = \begin{bmatrix}
\frac{E}{a(k+1)}-1 & -1\\
1 & 0
\end{bmatrix}
\quad V(k) = \begin{bmatrix} \frac{-b(k)}{a(k+1)} & 1-\frac{a(k)}{a(k+1)}\\
0 & 0
\end{bmatrix}.$$
Let
\begin{gather*}
\lambda_\pm (k) = \frac{E \pm \sqrt{E^2 - 4 a(k+1)^2}}{2
a(k+1)}\quad\quad \Lambda(k) = \begin{bmatrix}
\lambda_+&0\\0&\lambda_-\end{bmatrix} \\
 S(k) =\frac{1}{\lambda_+ -
\lambda_-}\begin{bmatrix}1&-\lambda_-\\-1&\lambda_+\end{bmatrix}.
\end{gather*}
Then $A(k) = S(k)^{-1}\Lambda(k)S(k)$.  If $|E|>2$ and $k$ is
large enough, then $\Lambda(k)$ satisfies the dichotomy condition
\eqref{dichotomy condition off}.  It is easy to see that the rest
of the hypotheses in Proposition \ref{lemma discrete HW off final}
are satisfied for $s=1$, so there are solutions of the form
$$\Psi_\pm (k) = S(k)^{-1}\Biggl( \prod_{l=k_0}^{k-1} \bigl(
\lambda_\pm (k) + \widetilde V(l)_\pm \bigr) \Biggr) \bigl( e_\pm
+ l^2_1\bigr)$$ where $\widetilde V_+ = \widetilde V_{11}$,
$\widetilde V_- = \widetilde V_{22}$, $e_+ = e_1$, $e_- = e_2$,
and $\|\widetilde V(l)\| \in l^2_1$.  We also have
\begin{equation*}
\lambda_\pm(k) + \widetilde V(k)_\pm = \lambda_\pm(k) \bigl(1 \pm
r_\pm (k) \bigr)
\end{equation*}
where $$r_\pm (k) = \frac{a(k+1) \bigl( \lambda_\pm(k) -
\lambda_\pm(k+1)\bigr) + \lambda_\mp(k) \bigl(a(k+1) - a(k) \bigr)
- b(k)}{a(k+1) \bigl(\lambda_\pm(k+1) - \lambda_\mp
(k+1)\bigr)}.$$  We now claim that
$$\sum_{l=k}^\infty r_\pm (l) \in l^2_1,$$ so in particular we can subsume the $1 +
r$ terms into the error to write $$\Psi_\pm (k) = c_\pm \Biggl(
\prod_{l=k_0}^{k-1} \lambda_\pm (l) \Biggr)\Bigl( S(k)^{-1}e_\pm +
l^2_1 \Bigr).$$

To see this is indeed the case, we make the following
observations.  First, $a(k) \rightarrow 1$, $\lambda(J), \kappa(J)
\in l^2_1$, and $\lambda_+ (k)$ and $\lambda_- (k)$ tend to
different finite constants.  In this way we see
$$\sum_{l=k}^\infty \Biggl(\frac{\lambda_\mp(l) \bigl(a(l+1) - a(l) \bigr)
- b(l)}{a(l+1) \bigl(\lambda_\pm(l+1) - \lambda_\mp
(l+1)\bigr)}\Biggr) \in l^2_1.$$  Second, we can write
$\lambda_+(k) - \lambda_+(k+1)$ as
\begin{align*}
\frac{\sqrt{E^2 - 4a_{k+1}^2}}{2a_{k+1}} &- \frac{\sqrt{E^2 -
4a_{k+2}^2}}{2a_{k+2}} = \frac{a_{k+2}\sqrt{E^2 -
4a_{k+1}^2}-a_{k+1}\sqrt{E^2 - 4a_{k+2}^2}}{2a_{k+1}a_{k+2}}\\
&=\sqrt{E^2 -
4a_{k+1}^2}\Biggl(\frac{a_{k+1}-a_{k+2}}{2a_{k+1}a_{k+2}}\Biggr)\\
&\quad+ a_{k+1}\frac{\sqrt{E^2 - 4a_{k+2}^2}-\sqrt{E^2 -
4a_{k+1}}}{a_{k+1}a_{k+2}}.
\end{align*}
Because $\kappa(J) \in l^2_1$, the first term is summable to be in
$l^2_1$ as well. To see the same is true of the second term, we do
a Taylor expansion of $\sqrt{E^2 - 4a^2}$ around the point
$E^2-4$.  After cancelling the constant terms we see that because
$\kappa(J) \in l^2_1$ we have
$$\sum_{k=n}^\infty a_{k+1}\Biggl(\frac{\sqrt{E^2 -
4a_{k+2}^2}-\sqrt{E^2 - 4a_{k+1}}}{a_{k+1}a_{k+2}}\Biggr) \in
l^2_1.$$  So the second term sums to be in $l^2_1$ as well,
proving the claim

Now, $E = \beta + \beta^{-1}$ and $$\beta^{\pm 1} = \frac{E \pm
\sqrt{E^2 - 4}}{2}$$ so
\begin{align*}
\lambda_\pm (k) &= \frac{\beta^{\pm 1}}{a(k+1)} \Biggl(\frac{E \pm
\sqrt{E^2 - 4a(k+1)^2}}{E \pm
\sqrt{E^2-4}}\Biggr)\\
&= \frac{\beta^{\pm 1}}{a(k+1)} \bigl( 1 \pm q_\pm(k)\bigr).
\end{align*}
Arguing as we did for the $r_\pm$ terms we find
$$\sum_{l=k}^\infty q_\pm (l) \in l^2_1,$$ so we can
subsume these products into the error term as well.  Finally,
using that $\kappa(J) \in l^2_1$ and taking the top row of
$\Psi_\pm$ we see
$$\psi_\pm(k) = c_\pm \beta^{\pm k} \bigl( 1 + l^2_1 \bigr),$$ as
claimed.
\end{proof}

\begin{proof}[Proof of Proposition \ref{lemma remove eits}(1)]
Recall that
\begin{gather*}
|\tilde \kappa(n-1) - \kappa(n-1)| = \Biggl| \sum_{k=n}^\infty
 a(k)^2 \Biggl(\frac{c(k-1)c(k+1)}{c(k)^2}-1\Biggr)\Biggr|\\
|\tilde \lambda(n-1) - \lambda(n-1) | =  |\gamma| \Biggl|
\sum_{k=n}^\infty
\Biggl(\frac{a(k-1)\phi(k-1)\phi(k)}{c(k-1)}-\frac{a(k) \phi(k)
\phi(k+1)}{c(k)}\Biggr)\Biggr|
\end{gather*}
where $J \phi = E \phi$, $\phi(0) = 0$, $\phi(1) = 1$ and $$c(n) =
1 + \gamma \sum_{j=1}^n |\phi(j)|^2.$$

Write $\phi$ as a linear combination of $\psi_+$ and $\psi_-$. Let
us first suppose that $\phi$ is just a multiple of $\psi_-$. As
$\psi_-$ is geometrically decreasing, the same is true of
$$a(k)^2 \Biggl(\frac{c(k-1)c(k+1)}{c(k)^2}-1\Biggr)$$ and
$$\Biggl( \frac{a(k-1)\phi(k-1)\phi(k)}{c(k-1)}-\frac{a(k) \phi(k)
\phi(k+1)}{c(k)}\Biggr).$$  So in this case it is easy to see that
$|\tilde \kappa(n-1) - \kappa(n-1)|$ and $|\tilde \lambda(n-1) -
\lambda(n-1) |$ are in $l^2_1$.

Now suppose that $\phi$ is not just a multiple of $\psi_-$. As
$\psi_+$ increases geometrically and $\psi_-$ decays
geometrically, we see
\begin{equation}\label{c equation}
c(k) \sim 1 + \gamma \sum_{l=1}^k \psi(l)^2 \sim 1 + \gamma
\sum_{l=1}^k \beta^{2l} \bigl( 1 + \tilde \delta(l) \bigr) \sim 1
+ \beta^{2k} \bigl( 1 + \delta(k) \bigr)
\end{equation}
where $\tilde \delta(k), \delta(k)$ represent some sequences in
$l^2_1$, and ``$\sim$'' indicates asymptotic equivalence (modulo
multiplication by constants). Similarly
$$\psi(k)\psi(k+1) \sim \beta^{2k+1}\bigl( 1 + \eps(k) \bigr)$$
for some $\eps \in l^2_1$. Combining these shows
\begin{align*}
\Biggl| &\frac{a(k-1)\psi(k-1)\psi(k)}{c(k-1)} -
\frac{a(k)\psi(k)\psi(k+1)}{c(k)}\Biggr|\\
&\lesssim \Biggl|\frac{a(k-1)\beta^{2k-1}\bigl(1 + \eps(k-1)
\bigr) -
a(k)\beta^{2k+1}\bigl(1 + \eps(k)\bigr)}{c(k-1)c(k)}\Biggr|\\
 &\quad+ \Biggl|\frac{\beta^{4k-1}\Bigl( a(k-1) \bigl(1 + \eps(k-1)
\bigr) - a(k)\bigl( 1 + \eps(k)\bigr) \Bigr)}{c(k-1)c(k)}\Biggr|\\
 &\quad+ \Biggl|\frac{\beta^{4k-1}\Bigl( a(k-1) \delta(k)\bigl(1 +
\eps(k-1) \bigr) - a(k)\delta(k-1)\bigl( 1 + \eps(k)\bigr)
\Bigr)}{c(k-1)c(k)}\Biggr|.
\end{align*}

Because $c(k-1)c(k) \sim \beta^{4k-1}$, the first term is
geometrically decreasing, so okay by Lemma \ref{lemma sum
products}.

Terms of the form
$$\frac{a(k-1)\beta^{4k-1}}{c(k-1)c(k)}\eps(k-1)\delta(k)$$ are in
$l^1_1$, being products of $l^2_1$ sequences.  Again, Lemma
\ref{lemma sum products} shows this is fine.

This leaves terms of the form
$$\frac{\beta^{4k-1}}{c(k-1)c(k)}
\bigl(\eps(k-1) - \eps(k)\bigr)$$ for some sequence $\eps \in
l^2_1$.  So it is sufficient to prove $$\sum_{k=n}^\infty \Biggl(
\frac{\beta^{4k-1}}{c(k-1)c(k)} \bigl(\eps(k-1) -
\eps(k)\bigr)\Biggr) \in l^2_1.$$

Let $$C(k) = \frac{\beta^{4k-1}}{c(k-1)c(k)}.$$  Summing by parts
shows
\begin{equation}\label{sum by parts}
\sum_{k=n}^\infty C(k) \bigl( \eps(k-1) - \eps(k) \bigr) =
C(n)\eps(n-1) + \sum_{k=n}^\infty \eps(k) \bigl( C(k+1) - C(k)
\bigr).
\end{equation}
The first term is clearly in $l^2_1$, so consider the second.
Using \eqref{c equation} we can write
\begin{align*}
|C(k+1) - &C(k)|\\ &= \Biggl| \frac{\beta^{4k+3}}{c(k)c(k+1)} -
\frac{\beta^{4k-1}}{c(k)c(k+1)}\Biggr|\\
&\lesssim \Biggl|\frac{\beta^{4k-1}}{c(k-1)c(k)c(k+1)} \Bigl(
\bigl( \beta^4-1\bigr) + \beta^{2k+2} \bigl( \delta(k-1) - \delta
(k+1)\bigr) \Bigr)\Biggr|.
\end{align*}
As $c(k-1)c(k)c(k+1) \sim \beta^{6k}$ the first term is
geometrically decaying and the second term is in $l^2_1$.
Combining this with \eqref{sum by parts} and Lemma \ref{lemma sum
products} shows that
$$\sum_{k=n}^\infty \Biggl( \frac{\beta^{4k-1}}{c(k-1)c(k)}
\bigl(\eps(k-1) - \eps(k)\bigr)\Biggr) \in l^2_1.$$

This completes the proof for the $\lambda$'s.  The proof for the
$\kappa$'s is similar and simpler.
\end{proof}

%
%
%
%

\section{Proof of Theorem \ref{main theorem} ($(2) \Rightarrow (1)$)}
By assumption, $J$ has finitely many eigenvalues, and they all lie
in $\R \setminus \interval$.  By Theorem \ref{lemma double
commutation} and Proposition \ref{lemma remove eits} we see it
suffices to prove the theorem when $\sigma (J) \subseteq
\interval$, which we now assume.

Now, $d\nu(x) = \chi_{[-2,2]}(x) v(x) dx$ and $v(x)$ has one of
the forms
$$\Bigl(\sqrt{4-x^2}\Bigr)^{\pm 1} v_0 (x) \quad \text{or}\quad
\Biggl(\sqrt{\frac{2-x}{2+x}}\Biggr)^{\pm 1} v_0 (x)$$ with $\log
v_0 \in \hoh$.  Define $$w(\theta) = c v_0(2 \cos \theta)$$ with
$c$ chosen to normalize $w$ to be a probability measure on
$\partial \D$. Notice that $$d\nu = Sz^\any (d\mu)$$ where $d\mu =
w \frac{d\theta}{2\pi}$ and $\any$ is one of $\even, \odd, \plus,
\minus$ according to which of the above forms $v(x)$ has.
Therefore, $\alpha$ and $a,b$ are related by one of the Direct
Geronimus Relations of Proposition \ref{lemma Geronimus
relations}.

By assumption $\log w \in \hoh$, so Theorem \ref{lemma
Golinskii-Ibragimov} shows that its Verblunsky coefficients
satisfy $\alpha \in l^2_1$.  By Proposition \ref{lemma Geronimus
relations} we see that
\begin{equation*}
\begin{split}
\kappa_n = \alpha_{2n-1} + K(\alpha)_n\\
\lambda_n = \alpha_{2n-2} + L(\alpha)_n
\end{split}
\end{equation*}
where $K(\alpha)_n$ and $L(\alpha)_n$ are sums from $n$ to
infinity of terms that are at least quadratic in $\alpha$. So by
Lemma \ref{lemma sum products} we see that $\lambda, \kappa \in
l^2_1$ too.

%
%
%
%

\section{Asymptotic integration redux}
For this section we will make the standing assumption that the
parameters defined by \eqref{lambda kappa def} obey
$$\kappa(J),\lambda(K) \in l^2_1.$$  In Section 7 we will need asymptotics on solutions at
energies $E = \pm2$. As before we will use asymptotic integration,
but because the recurrence equation at $E = \pm2$ yields a system
with a Jordan anomaly, we cannot use the results of Section 2.
Instead we construct a small solution $\psi_s$ and big solution
$\psi_b$:
\begin{proposition}\label{lemma edge spectrum asymptotics}
There are solutions $\psi_s$ and $\psi_b$ to $J \psi = E \psi$ at
energy $E = \pm 2$ such that
$$\Biggl| \frac{\psi_s(k)}{\psi_b(k)} \Biggr| \rightarrow 0$$
and
$$\pm \frac{\psi_s (k+1)}{\psi_s(k)} = 1 + l^2_1.$$  Moreover, for
either solution and for $k$ sufficiently large, $$(\pm1)^k \psi
(k)
> 0.$$
\end{proposition}

The rest of this section is devoted to a proof of this statement
for $E=2$, the proof for $E = -2$ being analogous. Recall we can
write the recurrence equation as
$$\Psi(k+1) = \begin{bmatrix}
\frac{2-b(k+1)}{a(k+1)} & -\frac{a(k)}{a(k+1)}\\
1 & 0\end{bmatrix} \Psi(k)$$ where $$\Psi(k) =
\begin{bmatrix}
\psi(k)\\
\psi(k-1)
\end{bmatrix}.$$

We begin with some preliminary transformations. Let $$S =
\begin{bmatrix} 1 & 1/2\\1&-1/2\end{bmatrix}$$ and
let $\Psi(k) = S \Phi(k)$. Then $\Phi$ solves $$\Phi(k+1) = [J +
B(k)]\Phi(k)$$ where
\begin{gather*}
J = \begin{bmatrix}1&1\\0&1\end{bmatrix}\\
B(k) =
\begin{bmatrix}\frac{\bigl(a(k+1)-1\bigr) + \bigl(a(k)-1\bigr) -
b(k)}{2a(k+1)} & \frac{-3\bigl(a(k+1)-1\bigr) + \bigl(a(k)-1\bigr)
- b(k)}{4 a(k+1)} \\ \frac{2\bigl(a(k+1)-1\bigr) +
\bigl(a(k)-1\bigr) - b(k)}{a(k+1)} & \frac{-3\bigl(a(k+1)-1\bigr)
+ \bigl(a(k)-1\bigr) - b(k)}{2a(k+1)}\end{bmatrix}.
\end{gather*}  In particular, notice that $$\Biggl\|\sum_{l=k}^\infty
B(l)\Biggr\| \in l^2_1.$$

\begin{lemma}\label{lemma jordan anomaly HL}
There exists a sequence of matrices $$Q(k) = \begin{bmatrix}
1&1\\q(k)&1\end{bmatrix}$$ such that $q \in l^2_1$,
\begin{equation}\label{HL equality jordan}
Q(k+1)^{-1}[J + B(k)]Q(k) = L(k) + M(k),
\end{equation}
and
\begin{equation*}
\|M(k)\| \in l^1_1,\quad L(k) = \begin{bmatrix}1 + \alpha(k) & 1 + \beta(k)\\
0&1 + \gamma(k) \end{bmatrix},
\end{equation*}
where
\begin{equation*}
\sum_{l=k}^\infty \alpha(l) \in l^2_1,\quad \sum_{l=k}^\infty
\beta(l) \in l^2_1,\quad \gamma(k) \in l^2_1.
\end{equation*}
In particular, if $\Phi(k) = Q(k)x(k)$ then $$x(k+1) = [L(k) +
M(k)]x(k).$$
\end{lemma}

\begin{proof}
Define $$Q(k) = \begin{bmatrix}1&1\\q(k)&1\end{bmatrix}$$ where
$q(k) = -\sum_{l=k}^\infty B(l)_{21} \in l^2_1$ by Lemma
\ref{lemma sum products}.  All the claimed properties are
straightforward calculations.
\end{proof}

As we seek asymptotics as $k \rightarrow \infty$, we need only
consider systems for $k$ larger than some $k_0$.  In particular,
we can choose $k_0$ so that $|\alpha(k)|,|\beta(k)|,|\gamma(k)|<
1$ for $k \geq k_0$.  In this case define
$$x(k) = P(k) z(k)$$ where $$P(k) =
\begin{bmatrix}1&0\\0&\prod_{j=1}^{k-1} \bigl( 1 + \gamma(j)
\bigr)\end{bmatrix}.$$ This transforms the $x$-system into
\begin{equation}\label{z equation}
z(k+1) = [J(k) + R(k)]z(k)
\end{equation}
where $$J(k) =
\begin{bmatrix}1 + \alpha(k) & 1 + \beta(k)\\0&1\end{bmatrix}$$ and $\|R(k)\| \in l^1_1$.
We will compare this to the simpler system $$y(k+1) = J(k)y(k).$$
We begin by finding a basis of solutions to the $y$-system.

\begin{lemma}\label{lemma jordan anomaly J basis}
The $y$-system above has two solutions $$y_s(k) =
\begin{bmatrix}u(k)\\0\end{bmatrix}\quad\text{and}\quad y_b(k) =
\begin{bmatrix} v(k)\\1\end{bmatrix}$$ such that
\begin{equation*}
u(k) = \prod_{j=1}^{k-1} \Bigl( 1 + \alpha(j) \Bigr)
\end{equation*}
\begin{equation}\label{u bdd below}
0 < |u(k)| \lesssim 1
\end{equation}
\begin{equation*}
 |v(k)|\sim k.
\end{equation*}
\end{lemma}

\begin{proof}
Let $u(1) = 1$ and $v(1) = 0$.  By the form of $J(k)$ we see that
$$u(k) = \prod_{j=1}^{k-1} \Bigl( 1 + \alpha(j) \Bigr).$$  Because
the $\alpha$'s are conditionally summable, the product defining
$u(k)$ converges to some finite number as $k \rightarrow \infty$,
so $|u(k)| \lesssim 1$. As we have assumed that $|\alpha(k)|<1$,
we also have $u(k) \neq 0$.

Now,
\begin{align}\label{v equation}
\begin{split}
v(k+1) &= \Bigl( 1 + \alpha(k) \Bigr) v(k) +
\Bigl( 1 + \beta(k) \Bigr)\\
& = \sum_{j=1}^k \frac{u(k+1)}{u(j+1)} \Bigl( 1 + \beta(j) \Bigr).
\end{split}
\end{align}
By \eqref{u bdd below} and $\beta(j) \rightarrow 0$ we have
$$|v(k)| \lesssim \sum_{l=1}^k 1 \lesssim k.$$  Moreover, there is some
$j_0$ so that for $k \geq j \geq j_0$, $$\frac{u(k+1)}{u(j+1)}
\Bigl( 1 + \beta(j) \Bigr)$$ is sign-definite. Without loss,
assume that it is positive, so for $j$ and $k$ large enough we
have
$$\frac{u(k+1)}{u(j+1)} \Bigl(1 + \beta(j)\Bigr)
 \gtrsim 1.$$  Thus, $|v(k)| \gtrsim k$ too.
\end{proof}

Now let $$Y(k) = \begin{bmatrix}u(k) & v(k)\\0&1\end{bmatrix}$$ be
a fundamental matrix for the $y$-system.  The next two lemmas
construct the desired solutions to \eqref{z equation}.

\begin{lemma}\label{lemma jordan anomaly small solution}
There is a bounded solution to the system \eqref{z equation} that
has
$$\|z(k+1) - z(k)\| \in l^2_1.$$ Moreover, $z(k)$ is sign-definite
for large enough $k$, and $\|z(k)\| > 0$.
\end{lemma}

\begin{proof}
Consider the operator
\begin{equation}\label{T jordan}
[Tz](k) = - \sum_{l=k}^\infty Y(k)Y(l+1)^{-1}R(l)z(l)
\end{equation}
 acting on
$l^\infty(\N;\C^2)$, with $k \geq k_1$ and $k_1$ to be chosen
momentarily. Notice that
$$Y(k)Y(l+1)^{-1} =
\begin{bmatrix}\frac{u(k)}{u(l+1)} & v(k) -
\frac{u(k)}{u(l+1)}v(l+1)\\0&1\end{bmatrix}.$$  By Lemma
\ref{lemma jordan anomaly J basis},
$$\Biggl| \frac{u(k)}{u(l+1)} \Biggr| \lesssim 1.$$ By \eqref{v
equation} we see that $$v(k) - \frac{u(k)}{u(l+1)}v(l+1) = -
\sum_{j=k}^l \frac{u(k)}{u(j+1)} \Bigl( 1 + \beta(j)\Bigr)$$ so
$$\Biggl| v(k) - \frac{u(k)}{u(l+1)}v(l+1) \Biggr| \lesssim
|l-k|.$$  Thus
\begin{equation}\label{YY bound}
\|Y(k)Y(l+1)^{-1}\| \lesssim |l-k|.
\end{equation}

If $z \in l^\infty$ we see $$\|[Tz](k)\|_\infty \lesssim
\|z\|_\infty\sum_{l=k_1}^\infty l \|R(l)\|.$$  Now, $\|R(l)\| \in
l^1_1$, so by choosing $k_1$ sufficiently large, we can ensure
$$\|Tz\| < \eps \|z\|$$ for some $\eps < 1$. Thus, $T$ is a
contraction on $l^\infty$, so in particular, given any $y \in
l^\infty$ there is a unique $z \in l^\infty$ solving $$z = y +
Tz.$$  If $y = y_s$ from Lemma \ref{lemma jordan anomaly J basis},
then by the form of $T$ (and a lengthy but easy calculation) we
see that this $z$ solves the $z$-equation.

Since $$\|Tz\| < \eps \|z\|$$ we see that $$\|y_s\| \leq \|z\| +
\|Tz\| < (1 + \eps)\|z\|$$ so $$\|z\| \gtrsim \|y_s\| > 0$$ by
Lemma \ref{lemma jordan anomaly J basis}.  Moreover, because
$y_s(k)$ is sign-definite for large enough $k$ and $\|Tz\| < \eps
\|z\|$, we see the same is true of $z$.

Next, notice that $$\|z(k+1) - z(k)\| \leq \|y_s(k+1) - y_s(k)\| +
\|[Tz](k+1) - [Tz](k)\|.$$  For the first term we use Lemma
\ref{lemma jordan anomaly J basis} to see
$$\|y_s(k+1) - y_s(k)\| \lesssim |\alpha(k)| \in
l^2_1.$$ For the second term we use \eqref{T jordan} to write
\begin{multline}\label{T jordan diff}
[Tz](k+1) - [Tz](k)\\ = Y(k)Y(k+1)^{-1}R(k)z(k) -
\sum_{l=k+1}^\infty [Y(k+1) - Y(k)]Y(l+1)^{-1}R(l)z(l).
\end{multline}

Now,
\begin{multline*}
[Y(k+1) - Y(k)]Y(l+1)^{-1} =\\
\begin{bmatrix} \frac{u(k+1) -
u(k)}{u(l+1)} & \Bigl( v(k+1) - \frac{u(k+1)}{u(l+1)}v(l+1) \Bigr)
+ \Bigl( v(k) - \frac{u(k)}{u(l+1)}v(l+1) \Bigr)\\
0&0\end{bmatrix}.
\end{multline*}
By Lemma \ref{lemma jordan anomaly J basis},
$$\Bigl| \frac{u(k+1) -
u(k)}{u(l+1)} \Bigr| \lesssim 1.$$  For the other term in the
matrix we use \eqref{v equation} to rewrite
$$\Bigl( v(k+1) - \frac{u(k+1)}{u(j+1)}v(l+1) \Bigr) + \Bigl( v(k) -
\frac{u(k)}{u(j+1)}v(j+1) \Bigr)$$
\begin{align*}
&= -\sum_{j=k+1}^l \frac{u(k+1)}{u(j+1)} \Bigl( 1 + \beta(j)\Bigr)
+ \sum_{j=k}^l \frac{u(k)}{u(j+1)} \Bigl( 1 + \beta(j) \Bigr)\\
&= \frac{u(k)}{u(k+1)} \Bigl( 1 + \beta(k)\Bigr) + \bigl( u(k) -
u(k+1) \bigr) \sum_{j=k+1}^l \frac{1}{u(j+1)} \Bigl( 1 +
\beta(j)\Bigr)\\
&= \frac{u(k)}{u(k+1)} \Bigl( 1 + \beta(k)\Bigr) - \alpha(k)
u(k)\sum_{j=k+1}^l \frac{1}{u(j+1)} \Bigl( 1 + \beta(j)\Bigr).
\end{align*}
In particular we have $\|[Y(k+1) - Y(k)]Y(l+1)^{-1}\| \lesssim 1 +
|\alpha(k)| l$.

Plugging this into \eqref{T jordan diff} and using $z \in
l^\infty$ and \eqref{YY bound} we find
\begin{align*}
\|[Tz](k+1) - [Tz](k)\| &\lesssim \|R(k)\| + \sum_{l=k+1}^\infty
\|R(l)\| + |\alpha(k)|\sum_{l=k+1}^\infty l\|R(l)\|\\
&\lesssim \|R(k)\| + \sum_{l=k+1}^\infty \|R(l)\| + |\alpha(k)|.
\end{align*}
The first and third terms are clearly $l^2_1$, and by Lemma
\ref{lemma sum products} so is the second.  Thus $\|z(k+1) -
z(k)\| \in l^2_1$.
\end{proof}

\begin{lemma}\label{lemma jordan anomaly big solution}
There is a solution $$z_b(k) =
\begin{bmatrix}z_{b1}(k)\\z_{b2}(k)\end{bmatrix}$$ to the $z$-system
that is sign-definite for $k$ large enough and has $|z_{b1}(k)|
\sim k$ and $|z_{b2}(k)|\lesssim 1$.
\end{lemma}

\begin{proof}
Again, we compare the $z$-system to the $y$-system and use Lemma
\ref{lemma jordan anomaly J basis}.  This time, consider the
operator
\begin{equation}\label{T jordan big}
[Tz](k) = \sum_{l=k_1}^{k-1} Y(k) Y(l+1)^{-1}R(l)z(l)
\end{equation}
with $k_1\geq 1$ to be chosen momentarily. Let $z_0 = y_b$ and
$z_{j+1} = y_b + Tz_j$. Then
$$\|z_{j+1}(k) - z_j(k)\| = \|T[z_j - z_{j-1}](k)\|$$ and
$$\|z_{1}(k) - z_0(k)\| = \|T[y](k)\|.$$  By \eqref{YY bound} and
Lemma \ref{lemma jordan anomaly J basis} we have
\begin{align*}
\|[Ty](k)\| &\leq \sum_{l=k_1}^{k-1}
\|Y(k)Y(l+1)^{-1}\|\cdot\|R(l)\|\cdot \|y(l)\|\\
&\lesssim k \sum_{l=k_1}^{k-1} \|R(l)\| l.
\end{align*}
We can choose $k_1$ sufficiently large that $$\|[Ty](k)\| < k
\eps$$ where $\eps < 1$.  Then inductively we find that
$$\|z_{j+1}(k) - z_j(k)\| < k \eps^{j+1}.$$  In particular, for
each $k$, $z_j(k) \rightarrow z_b(k)$ as $j \rightarrow \infty$
and $z_b = y_b + Tz_b$.

By the form of $T$ we see that because $y_b$ solves the
$y$-equation, $z_b$ solves the $z$-equation.  Moreover, as
$$\|[Tz_b](k)\| < k \eps$$ and $$\|y_b(k)\| \sim k$$ we have
\begin{equation}\label{z is big}
\|z_b\| \sim k.
\end{equation}
Finally, because $y_b(k)$ is sign-definite for large enough $k$
and $\|[Tz_b](k)\| < k \eps$, we see the same is true for $z_{b}$.

To deduce the component bounds, we expand
$$Y(k)Y(l+1)^{-1}R(l)z_b(l)$$ and notice that the bottom component
is bounded by $$|z_{b1}(l)R(l)_{21}| + |z_{b2}(l)R(l)_{22}|.$$
Plugging this into \eqref{T jordan big} shows
\begin{align*}
|z_{b2}(k)| &\lesssim \sum_{l=k_1}^{k-1} |z_{b1}(l)R(l)_{21}| + |z_{b2}(l)R(l)_{22}|\\
&\lesssim \sum_{l=k_1}^{k-1} l \|R(l)\| \lesssim 1.
\end{align*}
Combining this with \eqref{z is big} yields the final bound.
\end{proof}

\begin{proof}[Proof of Proposition \ref{lemma edge spectrum asymptotics}]
Undoing the transformations we find that $$\Psi(k) = S Q(k) P(k)
z(k)$$ and therefore that $$\psi(k) = \frac{1}{2}\Bigl( \bigl( 1 +
q(k)\bigr) z_1(k) + 2 \prod_{j=1}^{k-1} \bigl( 1 + \gamma(j)
\bigr) z_2(k) \Bigr)$$ where $$z(k) =
\begin{bmatrix}z_1(k)\\z_2(k)\end{bmatrix}.$$  Let $\psi_s$ and
$\psi_b$ correspond to taking $z$ to be $z_s$ and $z_b$.  All the
claimed properties now follow from Lemmas \ref{lemma jordan
anomaly small solution} and \ref{lemma jordan anomaly big
solution}.
\end{proof}

%
%
%
%

\section{The Inverse Geronimus Relations}
Recall that the Direct Geronimus Relations provide formulas for
$a^\any,b^\any$ in terms of $\alpha$.  In this section we go the
other way. We begin by determining whether a particular Jacobi
matrix is in the range of the \Szego maps based on the values of
its $m$-function. Note that while $Sz^\even$ maps onto all
probability measures supported on $\interval$, the ranges of the
other three maps are given by
\begin{gather}\label{szego map range}
\begin{split}
Ran(Sz^\odd) = \left\{ d\nu : \int_{-2}^2 \frac{ d\nu (x)}{4 -
x^2}
< \infty \right\}\\
Ran(Sz^\plusminus) = \left\{ d\nu : \int_{-2}^2 \frac{ d\nu (x)}{2
\mp x} < \infty \right\}.
\end{split}
\end{gather}

If $x \in \R$, write $$m(x + i0) = \lim_{\eps \downarrow 0} m(x +
i \eps),$$ and write $m(x)$ to indicate the value of the integral
$$\int\frac{d\nu(x)}{x-z}$$ (which may be infinite).

We begin by developing some elementary properties of the
$m$-functions, which we then use to study the associated
polynomials.

\begin{lemma}\label{lemma szego map range and m}
Let $J$ be a Jacobi matrix with $\sigma(J) \subseteq \interval$.
Then
\begin{align*}
J \in Ran(Sz^\odd) &\Leftrightarrow m(-2)-m(2) < \infty\\
J \in Ran(Sz^\plusminus) &\Leftrightarrow \mp m(\pm2) < \infty.
\end{align*}
\end{lemma}

\begin{proof}
The second line follows from the definition of the $m$-function
and \eqref{szego map range}. For the first line note
\begin{equation*}
m(-2) - m(2) = \int_{-2}^2
\Bigl(\frac{1}{2+x}+\frac{1}{2-x}\Bigr)d\nu(x) = 4 \int_{-2}^2
\frac{d\nu(x)}{4-x^2}
\end{equation*}
and again use \eqref{szego map range}.
\end{proof}

As with the ranges, the normalization constants \eqref{c
definition} have interpretations in terms of the $m$-function:

\begin{lemma}\label{lemma m alpha}
If $m^\any(x)$ is the $m$-function for $d\nu^\any$ then
\begin{gather*}
\mp m^\odd (\pm 2) = \frac{1}{(1 \mp \alpha_0)(1 - \alpha_1)}\\
\mp m^\plusminus(\pm2) = \frac{1}{2(1 \mp \alpha_0)}.
\end{gather*}
\end{lemma}

\begin{proof}
By \eqref{weight relation} we can write
\begin{gather*}
d\nu^\odd (x) = \frac{(2-x)(2+x)}{2(1-\alpha_0^2)(1 -
\alpha_1)}d\nu^\even (x)\\
d\nu^\plusminus (x) = \frac{2 \mp x}{2(1 \mp \alpha_0)}d\nu^\even
(x).
\end{gather*}
The values of $m^\plusminus$ then follow from $d\nu^\even$ being a
probability measure. For the $m^\odd$ values we have
\begin{align*}
m^\odd (-2) &= \int_{-2}^2 \frac{d\nu^\odd(x)}{2+x} =
\frac{1}{2(1-\alpha_0^2)(1-\alpha_1)} \int_{-2}^2 (2-x)d\nu^\even
(x)\\
&=\frac{1}{2(1-\alpha_0^2)(1-\alpha_1)} \int_0^{2\pi} 2 - (z +
z^{-1}) d\mu(z)\\
&= \frac{1}{(1-\alpha_0^2)(1-\alpha_1)}(1 -
\int_0^{2\pi}z d\mu(z))\\
&=\frac{1-\alpha_0}{2(1-\alpha_0^2)(1-\alpha_1)}.\\
\end{align*}
The value of $-m^\odd(2)$ follows similarly.
\end{proof}

We'll need lower bounds on the $m$-function:

\begin{lemma}\label{lemma m lower bound}
If $\sigma(J) \subseteq \interval$, then $\mp m(\pm 2) > 1/4.$
\end{lemma}

\begin{proof}
As $J$ has no eigenvalues off $[-2,2]$, $$m(E) = \int_{-2}^2
\frac{d\nu(x)}{x-E}.$$ For $t \in [-2,2]$ and $E > 2$, $t-E \geq
-4$. Because $d\nu$ is a probability measure that is not a point
mass at $t=2$, the Monotone Convergence Theorem implies
$$-m(2) = \lim_{E \downarrow 2}-m(E)> 1/4.$$
Similar arguments show $m(-2) > 1/4$.
\end{proof}

We now turn to the polynomials.  Given $d\mu$, write $P^\any_n(x)$
for the monic polynomial of degree $n$ with respect to the measure
$d\nu^\any = Sz^\any(d\mu)$. Similarly, let $Q^\any_n (x)$ be the
second-kind polynomial for $d\nu^\any$. That is, $Q$ solves the
same recurrence equation as $P$ but with initial conditions
$Q_{-1} \equiv -1$ and $Q_0 \equiv 0$. If $|m(x)|<\infty$, let
$F^\any_n (x) = m(x)P^\any_n(x) + Q^\any_n (x)$.

\begin{proposition}\label{lemma P F alpha}
Let $d\mu$ a nontrivial probability measure on $\partial \D$ that
is invariant under conjugation, and let $\alpha$ be its Verblunsky
parameters.  Then
\begin{gather*}
P^\even_{n+1}(2) = (1-\alpha_{2n-1})(1-\alpha_{2n})P^\even_n(2)\\
P^\even_{n+1}(-2) = -(1-\alpha_{2n-1})(1+\alpha_{2n})P^\even_n(-2)
\end{gather*}
\begin{gather*}
F^\odd_{n+1}(2) = (1+\alpha_{2n+1})(1+\alpha_{2n+2})F^\odd_n(2)\\
F^\odd_{n+1}(-2) = -(1+\alpha_{2n+1})(1-\alpha_{2n+2})F^\odd_n(-2)
\end{gather*}
\begin{gather*}
F^\plus_{n+1}(2) = (1+\alpha_{2n})(1+\alpha_{2n+1})F^\plus_n(2)\\
P^\plus_{n+1}(-2) = -(1+\alpha_{2n})(1-\alpha_{2n+1})P^\plus_n(-2)
\end{gather*}
\begin{gather*}
P^\minus_{n+1}(2) = (1-\alpha_{2n})(1-\alpha_{2n+1})P^\minus_n(2)\\
F^\minus_{n+1}(-2) =
-(1-\alpha_{2n})(1+\alpha_{2n+1})F^\minus_n(-2).
\end{gather*}
\end{proposition}

\begin{proof}
The proof is by induction.  As the arguments for any of the $P$'s
are virtually identical, we only present the proof for the case $P
= P^\even$.  Similarly, we only present the argument for the $F$'s
in the case $F = F^\plus$.

The desired relationship between $P_0 \equiv 1$ and $P_1(x) = x
-b_1$ follows from Proposition \ref{lemma Geronimus relations}:
\begin{equation*}
P_1(2) = 2-b_1 = 2 - \bigl(\alpha_{0}(1-\alpha_{-1}) \bigr) =
2(1-\alpha_0) =(1-\alpha_{-1})(1-\alpha_{0})P_0(2).
\end{equation*}

To deduce the desired relationship between $F_1(2)$ and $F_0(2)$,
we argue as follows.  By Lemma \ref{lemma m lower bound} and Lemma
\ref{lemma szego map range and m} we have that $\frac{1}{4} <
-m(2) < \infty$, so $F_0(2) = m(2)$. Next, recall that $P_{-1}
\equiv 0$, $P_{0} \equiv 1$, $Q_{-1} \equiv -1$, and $Q_{0} \equiv
0$. So
\begin{align*}
F_1(2) &= m(2)P_1(2) + Q_1(2)\\
&= -\frac{(2-b_1)}{2(1-\alpha_0)} + 1
=-\frac{2\alpha_0 - b_1}{2(1-\alpha_0)}\\
&=(1+\alpha_0)(1+\alpha_1)\frac{-1}{2(1-\alpha_0)} =
(1+\alpha_0)(1+\alpha_1)F_0(2)
\end{align*}
where we have used Proposition \ref{lemma Geronimus relations} and
Lemma \ref{lemma m alpha}.

Now, assume the formulas hold up to index $n-1$.  As $P_n$
satisfies the three-term recurrence equation we have
\begin{align*}
P_{n+1}(2) &= (2-b_{n+1})P_n(2) -
a_n^2 P_{n-1}(2)\\
&= \Bigl( (2-b_{n+1}) -
\frac{a_n^2}{(1-\alpha_{2n-3})(1-\alpha_{2n-2})}
\Bigr)P_n(2)\\
&=(1-\alpha_{2n-1})(1-\alpha_{2n})P_n(2)
\end{align*}
where the second equality is by the inductive hypothesis, and the
third equality is by Proposition \ref{lemma Geronimus relations}.

Similarly, $F_n$ satisfies the three-term recurrence equation, so
the same argument works:
\begin{align*}
F_{n+1}(2) &= (2-b_{n+1})F_n(2) -
a_n^2 F_{n-1}(2)\\
&= \Bigl( (2-b_{n+1}) -
\frac{a_n^2}{(1+\alpha_{2n-2})(1+\alpha_{2n-1})}
\Bigr)F_n(2)\\
&=(1+\alpha_{2n})(1+\alpha_{2n+1})F_n(2).
\end{align*}
\end{proof}

\begin{proposition}[Inverse Geronimus Relations\footnote{The
case $d\nu = Sz^\even(d\mu)$ is due to \cite{Geronimus} (with an
alternate proof given in \cite{Damanik Killip}).  The statement in
the other three cases appears to be new (although anticipated in
\cite{Simon} and related to some formulas of
\cite{spaniards}).}]\label{lemma Inverse Geronimus relations} Let
$d\mu$ a nontrivial probability measure on $\partial \D$ that is
invariant under conjugation, and let $\alpha$ be its Verblunsky
parameters. Define
\begin{gather*}
A^\any_n = -\frac{P^\any_{n+1}(-2)}{P^\any_n(-2)} \quad B^\any_n =
\frac{P^\any_{n+1}(2)}{P^\any_n(2)}\\
C^\any_n = -\frac{F^\any_{n+1}(-2)}{F^\any_n(-2)} \quad D^\any_n =
\frac{F^\any_{n+1}(2)}{F^\any_n(2)}.
\end{gather*}
If $d\nu = Sz^\even(d\mu)$
\begin{equation*}
\alpha_{2n} = \frac{A^\even_n - B^\even_n}{A^\even_n + B^\even_n}
\quad \alpha_{2n-1} = 1 - \frac{1}{2}(A^\even_n + B^\even_n).
\end{equation*}
If $d\nu = Sz^\odd(d\mu)$
\begin{equation*}
-\alpha_{2n+2} = \frac{C^\odd_n - D^\odd_n}{C^\odd_n + D^\odd_n}
\quad -\alpha_{2n+1} = 1 - \frac{1}{2}(C^\odd_n + D^\odd_n).
\end{equation*}
If $d\nu = Sz^\plus(d\mu)$
\begin{equation*}
-\alpha_{2n+1} = \frac{A^\plus_n - D^\plus_n}{A^\plus_n +
D^\plus_n} \quad -\alpha_{2n} = 1 - \frac{1}{2}(A^\plus_n +
D^\plus_n).
\end{equation*}
If $d\nu = Sz^\minus(d\mu)$
\begin{equation*}
\alpha_{2n+1} = \frac{C^\minus_n - B^\minus_n}{C^\minus_n +
B^\minus_n} \quad \alpha_{2n} = 1 - \frac{1}{2}(C^\minus_n +
B^\minus_n).
\end{equation*}
\end{proposition}

By Sturm oscillation theory and that $Sz^\any (d\mu)$ is supported
in $\interval$, we see $(\pm1)^{n+1}P^\any_n(\pm2)$ and
$-(\pm1)^{n+1}F^\any_n(\pm2)$ are strictly positive for all $n >
0$.  In particular, the above ratios are all defined.

\begin{proof}
This is a simple calculation based on Proposition \ref{lemma P F
alpha}.
\end{proof}

%
%
%
%

\section{Some Weyl theory}
By Proposition \ref{lemma Inverse Geronimus relations} we see that
decay of the Verblunsky parameters is controlled by decay of the
sequences $A_n$, $B_n$, $C_n$, and $D_n$.  By Proposition
\ref{lemma edge spectrum asymptotics} we see that there is a
solution at $E=\pm2$ with the desired asymptotics.  The following
result connects these two ideas.

\begin{proposition}\label{lemma mABCD}
Let $J$ be a Jacobi matrix with $\sigma (J) \subseteq \interval$,
and let $A_n$, $B_n$, $C_n$, $D_n$ be defined as above. Then
$m(-2) < \infty$ implies $C_n = 1 + l^2_1$, and $m(-2) = \infty$
implies $A_n = 1 + l^2_1$. Similarly, $-m(2) < \infty$ implies
$D_n = 1 + l^2_1$, and $-m(2) = \infty$ implies $B_n = 1 + l^2_1$.
\end{proposition}

Let us write $p_n$ and $q_n$ for the orthonormal versions of $P_n$
and $Q_n$, and then $f_n(z) = m(z) p_n(z) + q_n(z)$. Proposition
\ref{lemma mABCD} is a trivial consequence of

\begin{proposition}\label{lemma mFP}
Let $J$ be a Jacobi matrix with $\sigma (J) \subseteq \interval$.
Then $m(-2) < \infty$ implies $(-1)^{n+1}f_n (-2) = s + l^2_1$,
and $m(-2) = \infty$ implies $(-1)^{n+1}p_n (-2) = s + l^2_1$, for
some $s \in \R$. Similarly, $-m(2) < \infty$ implies $f_n (2) = s
+ l^2_1$, and $-m(2) = \infty$ implies $p_n (2) = s + l^2_1$.
\end{proposition}

To prove this, we will use some Weyl theory.  Recall $p_n(z)$ and
$q_n(z)$ are solutions to $Ju = zu$ with $p_{-1} = q_0 = 0$ and
$p_0 = -q_{-1} = 1$.  When $z \in \C \setminus \R$, the Weyl
solution $f_n(z) = m(z)p_n(z) + q_n(z)$ is defined and satisfies
\begin{equation}\label{weyl norm}
\|f_n(z)\|^2_{l^2} = \frac{\im m(z)}{\im z}.
\end{equation}

As the $m$-function and the solutions $p$ and $q$ will play
prominent roles, we develop some of their key properties.  To
start, we relate the values of $m$ at $\pm2$ to its values at $\pm
2 + i \eps$.

\begin{lemma}\label{lemma m(2) m(2+ieps)}
Let $J$ be a Jacobi matrix with $\sigma (J) \subseteq \interval$,
$m$-function m, and spectral measure $d\nu$.  Then
\begin{gather*}
\int_{-2}^2 \frac{d\nu(t)}{2+t} < \infty \quad\Rightarrow\quad
m(-2) = m(-2 + i0)\\
\int_{-2}^2 \frac{d\nu(t)}{2+t} = \infty \quad\Rightarrow\quad
|m(-2 + i0)| = \infty.
\end{gather*}
Similarly,
\begin{gather*}
\int_{-2}^2 \frac{d\nu(t)}{2-t} < \infty \quad\Rightarrow\quad
m(2) =
m(2 + i0).\\
\int_{-2}^2 \frac{d\nu(t)}{2-t} = \infty \quad\Rightarrow\quad
|m(2 + i0)| = \infty.
\end{gather*}
In particular, when $m(\pm 2)$ is finite, we may write $m(\pm2)$
for $m(\pm 2 + i0)$ and then $f_n(\pm 2)$ for $f_n(\pm 2 + i0)$.
\end{lemma}

Notice that by Lemma \ref{lemma m lower bound}, $\mp m(\pm2)$ can
only diverge to $+\infty$.

\begin{proof}
The first implication follows from the Dominated Convergence
Theorem applied to $$-m(2+i\eps) = \int_{-2}^2
\frac{v(t)dt}{(2-t)+i\eps}.$$ The second implication follows from
the Monotone Convergence Theorem applied to $$-\re m(2 + i\eps) =
\int_{-2}^2 \frac{2-t}{(2-t)^2+\eps^2} d\nu(t).$$
\end{proof}

If $L \in \N$ and $u(n;z)$ solves $Ju=zu$, we define
$$\|u(z)\|_L^2 = \sum_{l=0}^L |u(l;z)|^2.$$  For non-integer
values of $L$ we define $\|u(z)\|_L$ to be the linear
interpolation between the values at $\lfloor L \rfloor$ and
$\lceil L \rceil$. Now suppose $x \in \R$ is fixed and $m(x+i0)$
exists finitely. Let $\eps, y' > 0$ be related by $$\sup_{0 < y
\leq y'} |m(z) - m(x+i0)| + y' = \frac{\eps^2}{4}$$ where $z =
x+iy$.  Note that $y'$ is a monotone function of $\eps$ that goes
to zero as $\eps$ goes to zero. Define $L(\eps)$ by
$$|y'|^{1/2} \|p_n(z')\|_{L(\eps)} = 1$$ where $z' = x + i y'$.
For each $y' > 0$, $L(\eps)$ exists because $p_{n}(z')$ is not in
$l^2$.

The following lemma is a discrete analog of Lemma 9 of
\cite{Gilbert Pearson}.  The proof is a direct translation, so we
omit it.

\begin{lemma}\label{lemma um u1 ratio}
Let $x \in \R$ and suppose that $m(x+i0)$ exists finitely.  Then
$$\frac{\|f_n(x+i0)\|_{L(\eps)}}{\|p_n(x)\|_{L(\eps)}}<\eps$$
whenever $\eps$ is sufficiently small.
\end{lemma}

Next we recall a result of \cite{Jitomirskaya Last}.

\begin{lemma}\label{lemma jito last}
Let $x \in \R$ and define $\tilde L(\eps)$ by
$$\|p_n(x)\|_{\tilde L(\eps)}\|q_n(x)\|_{\tilde L(\eps)} =
\frac{1}{2\eps}.$$  Then $\tilde L(\eps)$ is a well defined,
monotonely decreasing continuous function that goes to infinity as
$\eps$ goes to 0, and $$\frac{5-\sqrt{24}}{|m(x+i\eps)|} \leq
\frac{\|p_n(x)\|_{\tilde L(\eps)}}{\|q_n(x)\|_{\tilde L(\eps)}}
\leq \frac{5+\sqrt{24}}{|m(x+i\eps)|}.$$
\end{lemma}

\begin{proof}[Proof of Proposition \ref{lemma mFP}]
Again, we will only prove the statements for $E = 2$.  Suppose
first that $1/4 \leq -m(2) < \infty$.  Then by Lemma \ref{lemma
m(2) m(2+ieps)}, $m(2 + i 0)$ is finite and nonzero too.  So by
Lemma \ref{lemma jito last} we have that
$$\frac{\|p_n(2)\|_L}{\|q_n(2)\|_L}$$ remains finite and nonzero
as $L \uparrow \infty$.  As solutions at $E=2$ are of the form
$c_1 \psi_b + c_2 \psi_s$ for some $c_i \in \R$, we see that
$p_n(2)$ and $q_n(2)$ must be simultaneously bounded or
simultaneously unbounded. Because $p_{n}(2)$ and $q_{n}(2)$ form a
basis for solutions at $E=2$, we see they cannot both be bounded.
So they are both unbounded.

Now, by Lemma \ref{lemma um u1 ratio} we see that
$$\frac{\|f_n(2)\|_L}{\|p_n(2)\|_L} \rightarrow 0$$ as $L \uparrow
\infty$.  So $f_{n}(2)$ cannot be unbounded, and so has the form
$c \psi_s$ for some $c \in \R$. Now Proposition \ref{lemma edge
spectrum asymptotics} yields the desired result.

Now suppose that $-m(2) = \infty$.  Then by Lemma \ref{lemma m(2)
m(2+ieps)} we have that $|m(2 + i0)| = \infty$ too.  Then by Lemma
\ref{lemma jito last} we have
$$\frac{\|p_n(2)\|_L}{\|q_n(2)\|_L} \rightarrow 0$$ so we must
have that $p_{n}(2)$ remains bounded.  Thus, $p_{n}(2) = c
\psi_s(n)$ for some $c \in \R$, so again we are done by
Proposition \ref{lemma edge spectrum asymptotics}.
\end{proof}

%
%
%
%

\section{Proof of Theorem \ref{main theorem} ($(1) \Rightarrow (2)$)}
By Proposition \ref{lemma edge spectrum asymptotics} we see that
all solutions at $E = \pm 2$ eventually satisfy $(\pm 1)^k \psi(k)
> 0$.  So by the Sturm oscillation theorem for Jacobi matrices
(see chapter 4 of \cite{Teschl}), $J$ has only finitely-many
eigenvalues, all lying in $\R \setminus \interval$. So b$\beta(j)
\rightarrow 0$ Proposition \ref{lemma remove eits} it suffices to
prove the theorem when $\sigma (J) \subseteq \interval$, which we
now assume.

Consider the values of the $m$ function at $E = \pm 2$. We have
four cases:

\underline{Case 1}: $m(-2) = -m(2) = \infty$. As $Sz^\even$ is
onto, $d\nu \in Ran(Sz^\even)$, so choose $$R_n(-2) = P_n(-2)
\quad\quad R_n(2) = P_n(2) \quad\quad d\mu =
[Sz^\even]^{-1}(d\nu).$$

\underline{Case 2}: $m(-2) , -m(2) < \infty$. By Lemma \ref{lemma
szego map range and m}, $d\nu \in Ran(Sz^\odd)$, so choose
$$R_n(-2) = F_n(-2) \quad\quad R_n(2) = F_n(2) \quad\quad d\mu =
[Sz^\odd]^{-1}(d\nu).$$

\underline{Case 3}: $m(-2) = \infty$, $-m(2) < \infty$. By Lemma
\ref{lemma szego map range and m}, $d\nu \in Ran(Sz^\plus)$, so
choose
$$R_n(-2) = P_n(-2) \quad\quad R_n(2) = F_n(2) \quad\quad d\mu = [Sz^\plus]^{-1}(d\nu).$$

\underline{Case 4}: $m(-2) < \infty$, $-m(2) = \infty$. By Lemma
\ref{lemma szego map range and m}, $d\nu \in Ran(Sz^\minus)$, so
choose
$$R_n(-2) = F_n(-2) \quad\quad R_n(2) = P_n(2) \quad \quad d\mu = [Sz^\minus]^{-1}(d\nu).$$

In any case, let $\alpha$ be the Verblunsky parameters associated
to $d\mu$.  By Proposition \ref{lemma mABCD} we see that
$$\frac{R_{n+1}(-2)}{R_n(-2)} = 1 + l^2_1 \quad\quad \frac{R_{n+1}(2)}{R_n(2)} = 1 +
l^2_1.$$  Then by Proposition \ref{lemma Inverse Geronimus
relations} we see that $\alpha \in l^2_1$.  By Theorem \ref{lemma
Golinskii-Ibragimov} we see $\log w \in \hoh$, so by \eqref{weight
relation} we see $v \in \W$.

%
%
%
%

\end{document}